\documentclass[12pt]{amsart}
\usepackage[english]{babel}
\usepackage{amsmath,amscd,amsthm,amsfonts,latexsym}

\usepackage{fullpage}

\theoremstyle{plain}
\newtheorem{theorem}                 {Theorem}      [section]

\theoremstyle{definition}
\newtheorem{example}      [theorem]  {Example}
\newtheorem{remark}       [theorem]  {Remark}
\newtheorem{definition}   [theorem]  {Definition}

\newcommand{\dd}{\mathrm d}
\def\lra{\longrightarrow}
\def \Re{\mbox{${\mathcal Re}$}}

\def \r{\mbox{${\mathbb R}$}}
\def \C{\mbox{${\mathbb C}$}}
\def \L {\mathbb{L}}
\def \K {\mathbb{K}}
\def \S {\mathbb{S}}

\def \H {\mathbb{H}}
\def \s {\mathfrak{s}}
\def \h {\mathfrak{h}}

\def \a {\mathfrak{a}}
\def \s {\mathfrak{s}}

\begin{document}


\title{Minimal surfaces in $4$-dimensional Lorentzian Damek-Ricci spaces}
\author{Adriana A. Cintra}
\address{Departamento de Matem\'{a}tica, C.P. 03\\ UFG, 75801-615, Jata\'{i}, GO \\Brasil }

\email{adriana.cintra@ufg.br}

\author{Francesco Mercuri}
\address{Departamento de Matem\'{a}tica, C.P. 6065\\
IMECC, UNICAMP, 13081-970, Campinas, SP\\ Brasil}

\email{mercuri@ime.unicamp.br}

\author{Irene I. Onnis}

\address{Departamento de Matem\'{a}tica, C.P. 668\\ ICMC,
USP, 13560-970, S\~{a}o Carlos, SP\\ Brasil}

\email{onnis@icmc.usp.br}

\thanks{Work partially supported by Capes and CNPq, Brasil.}

\begin{abstract}
In this paper we will construct a Weierstrass type representation for minimal surfaces in $4$-dimensional Lorentzian Damek-Ricci spaces and we give some examples of such surfaces.

\end{abstract}

\maketitle
\section{Introduction}

Damek-Ricci spaces are semidirect products of Heisenberg groups with the real line. They were considered in \cite{DR} (see also \cite{Tri}), equipped with a left-invariant Riemannian metric, to give a negative answer, in high dimensions, to the question posed by Lichnrowicz: ``is a harmonic
Riemannian manifold necessarily a symmetric space?''

Beside a left-invariant  Riemannian metric, these spaces may be equipped with left-invariant Lorentzian metrics in essentially two ways: a Riemannian metric on the Heisenberg factor and a negative metric on the $\r$ factor, or a Lorentzian metric in the Heisenberg factor and a positive metric on $\r$. The aim of this paper is to study a Weierstrass representation for simply connected minimal surfaces in these spaces, in dimension four.

In \cite{MMP} the authors give a Weierstrass representation theorem for minimal surfaces in Riemannian manifolds. In \cite{Liramm} this representation has been extended for timelike and spacelike minimal surfaces in 3-dimensional Lorentzian manifolds. The results can be easily extended to the case of minimal surfaces  in Lorentzian manifolds of higher dimension. Most of the applications and examples of these results are given for 3-dimensional ambient spaces. For higher dimension, there is an application of this formula for minimal surfaces in 4-dimensional Damek-Ricci spaces equipped with a left-invariant Riemannian metric (see \cite{kotod}).

This paper is organized as follows: in Section~2 we describe the geometry of the Damek-Ricci spaces and, then, in the next section we discuss the extensions of the Weierstrass representation theorem for minimal surfaces in Riemannian and Lorentzian manifolds. In Sections~4 and 5 we adapt the Weierstrass representation to our situation. Finally we give examples of spacelike and timelike minimal surfaces in these spaces (in the 4-dimensional case).

\section{The geometry of the Damek-Ricci spaces}

\subsection{The generalized Heisenberg group}

Let $b_m$ and $z_n$ be real vector spaces of dimensions $m$ and $n$
respectively, and $\beta: b_m \times b_m \rightarrow z_n$ a skew-symmetric bilinear map.
In the direct sum $\h_{m + n} = b_m\oplus z_n$ we define the bracket
\begin{equation}\label{eq:4.2}
    [U + X, V + Y] = \beta(U,V).
\end{equation}
This product defines a Lie algebra structure on $\h_{m+n}$, whose center contains $z_n$.

We endow $b_m$ with a positive inner product and $z_n$ with a positive or Lorentzian inner product. We will denote by  $\langle \cdot,\cdot\rangle_{\h_{m + n}}$ the product metric. For $Z \in z_n$, we define $J_Z \in End(b_m)$ by
\begin{equation}\label{eq:4.1}
  \langle J_Z U, V\rangle_{\h_{m + n}} = \langle \beta(U,V),Z\rangle_{\h_{m + n}},
\end{equation}
for all $U,V \in b_m$ and $Z\in z_n$.

The Lie algebra $\h_{m + n}$ is called a {\em generalized Riemannian Heisenberg algebra} if the inner product in $z_n$ is positive and
$$ J_Z^2 = -\langle Z,Z\rangle_{\h_{m + n}}id_{b_m},$$
for all $Z \in z_n$. The associated simply connected Lie group, with the left-invariant metric, is called a {\em generalized Riemannian Heisenberg group}.

The Lie algebra $\h_{m + n}$  is called a {\em generalized Lorentzian Heisenberg algebra} if the inner product in $z_n$\ is Lorentzian and
$$ J_Z^2 = -\langle Z,Z\rangle_{\h_{m + n}}id_{b_m},\qquad \text{if}\;Z\;  \text{is spacelike,}$$
$$J_Z^2 = \langle Z,Z\rangle_{\h_{m + n}}id_{b_m},\qquad \text{if}\;Z\;\text{is timelike}.$$
The associated simply connected Lie group, with the left-invariant metric, is called a {\em generalized Lorentzian Heisenberg group.}

\subsection{Lorentzian Damek-Ricci spaces of the first kind}

Take the direct sum $\s_{m+n+1} = \h_{m+n}\oplus \a$, where $\a$ is a Lorentzian one-dimensional space and $\h_{m+n}$ is a generalized Riemannian Heisenberg algebra. A vector in $\s_{m+n+1}$ can be written in a unique way as $U + X + s\,A$, for some
$U \in b_m$, $X \in z_n$, $s \in \r$ and  a  non zero fixed vector $A$ in $\a$.

Given $U + X + r\,A, V + Y + s\,A\in \s_{m + n + 1}$, we define
$$
    \langle U + X + r\,A, V + Y + s\,A\rangle = \langle U + X,V + Y \rangle_{\h_{m+n}} - \,r\,s
$$
and
$$
    [U + X  + r\,A, V + Y + s\,A] = [U,V]_{\h_{m+n}} + \frac{1}{2}r\,V -\frac{1}{2}s\,U +r\,Y -s\,X.
$$

Then $\langle \cdot,\cdot\rangle$ is a Lorentzian metric and $[\cdot,\cdot]$ is a Lie bracket in $\s_{m + n + 1}$. Therefore, $\s_{m + n + 1}$ is a Lie algebra. Moreover, $\langle A,A \rangle = -1$.

\begin{definition}
The simply connected Lie group associated to $\s_{m + n +1}$, endowed with the induced left-invariant Lorentzian metric, is called a {\em Lorentzian Damek-Ricci space of the first kind} and will be denoted by $\S_{m + n + 1}^1$.
\end{definition}
The Levi-Civita connection $\nabla$ of $\S_{m + n + 1}^1$ is given by
$$
\begin{aligned}
    &\nabla_{V + Y + s\,A} (U + X + r\,A)=\\
    & - \frac{1}{2}\{J_Y U + J_X V + rV + [U,V] + 2rY + \langle U, V\rangle A + 2\langle X, Y\rangle A\}.
    \end{aligned}
$$

\subsection{Lorentzian Damek-Ricci spaces of the second kind}

We consider again the direct sum $\s_{m+n+1} = \h_{m+n}\oplus \a$, where now $\a$ is a Riemannian 1-dimensional space and $\h_{m+n}$ is a generalized Lorentzian Heisenberg algebra. The bracket is given as above and the metric is given by
$$
    \langle U + X + r\,A, V + Y + s\,A\rangle = \langle U + X,V + Y \rangle_{\h_{m+n}} + \,r\,s.
$$

\begin{definition} The simply connected Lie group associated to $\s_{m + n +1}$, endowed with the induced left-invariant Lorentzian metric, is called a {\em Lorentzian Damek-Ricci space of the second kind} and will be denoted by $\S_{m + n + 1}^{m+n}$.
\end{definition}

The Levi-Civita connection $\nabla$ of $\S_{m + n + 1}^{m+n}$ is given by
$$
\begin{aligned}
    &\nabla_{V + Y + s\,A} (U + X + r\,A)=\\
    & - \frac{1}{2}\{J_Y U + J_X V + rV + [U,V] + 2rY - \langle U, V\rangle A - 2\langle X, Y\rangle A\}.
    \end{aligned}
$$

\section{The Weierstrass representation}

The Weierstrass representation theorem is an important tool in the study of minimal surfaces  in $\r^n$ since it allows to bring in the powerful theory of holomorphic functions. The local version has been extended to the case of minimal surfaces in a Riemannian manifold in \cite{MMP} and for Lorentzian manifolds in \cite{Liramm}. In this section we will  briefly discuss such extensions. We start with the Riemannian case. Since the considerations are local, we can suppose that the ambient manifold is $\r^n$ with a Riemannian   metric $g = [g_{ij}]$.
\begin{theorem}\label{classicalwr} Let $\Omega \subseteq \C$ be an open set and let $f: \Omega \lra \r^n$ be a conformal minimal immersion. Let $\{u, v\}$ be conformal coordinates in $\Omega$\ and $z = u + iv$. Consider the complex tangent vector
$$\frac{\partial f}{\partial z} := \frac{1}{2}(\frac{\partial f}{\partial u} + i\frac{\partial f}{\partial v}),$$
where $i = \sqrt{-1}$. Let
$$\frac{\partial f}{\partial z} = \sum_{i=1}^n \phi_i \frac{\partial}{\partial x_i}.$$
Then
\begin{enumerate}
\item $\sum_{i,j=1}^n g_{ij} \phi_i\,\phi_j \neq 0$,
\item $\sum_{i,j=1}^n g_{ij}\phi_i \overline{\phi_j} = 0$,
\item $\displaystyle{\frac{\partial \phi_i}{\partial \overline{z}} + \sum_{j,l=1}^n \Gamma^i_{jl}\ \phi_j \overline{\phi_l} = 0},$
\end{enumerate}
where $\Gamma^i_{jl}$ are the Christoffel symbols of $g$.
Moreover, if $\Omega$ is simply connected, the functions
$$f_i := 2\, \Re\int  \phi_i\, \dd z$$
are well defined and define a conformal minimal immersion with complex tangent vector
 $$\displaystyle{\frac{\partial f}{\partial z} = \sum_{i=1}^n \phi_i \frac{\partial}{\partial x_i}}.$$
\end{theorem}

\begin{remark}
In the Theorem~\ref{classicalwr} the first condition guarantees that $f$ is an immersion, the second one that $f$ is conformal and the last one that $f$ is minimal.
\end{remark}

\begin{remark}
The third condition is called the {\em harmonicity condition} since it just says that the tension field vanishes.
\end{remark}

In the case that $g$ is a Lorentzian metric the essential difference is the following: for spacelike surfaces (i.e. if $f^*g$ is Riemannian) the statement is the same. For timelike surfaces (i.e. if $f^*g$ is Lorentzian) the expression $\displaystyle{\frac{\partial \phi_i}{\partial \overline{z}}}$, as well as the conjugation, has to be understood in the Lorentz or {\em paracomplex} sense.

We recall that the algebra of {\em paracomplex numbers} is the algebra
$$\L = \{ a + \tau b,\ a, b \in \r\},$$
where $\tau$ is an imaginary unit with $\tau^2 = 1$. The operations are the obvious ones and the set of zero divisors is the set
$$K = \{ a \pm \tau a,\,a \in \r\}.$$
This algebra is isomorphic to $\r \oplus \r$ via the map $$a + \tau b \longmapsto \frac{1}{2}(a+b, a-b).$$
Paraconjugation and norm are defined as in the complex case and $z \in \L \setminus K$ is invertible with inverse $\displaystyle z^{-1} = \bar{z}/(z\bar{z})$.

The set $\L$ has a natural topology as a $2$-dimensional real vector space.
\begin{definition}
Let $\Omega\subseteq \L$ be an open set and $z_0 \in \Omega$.
The $\L$-derivative of a function $f:\Omega\rightarrow\L$ at $z_0$ is defined by
$$
f'(z_0):= \lim_{z\rightarrow z_{0} \atop{z - z_{0} \in \L\setminus K\cup\{0\}}}\frac{f(z) - f(z_0)}{z - z_0},
$$
if the limit exists. If $f'(z_0)$ exists, we will say that $f$ is $\L$-differentiable
at $z_0$.
\end{definition}

\begin{remark}
The condition of $\L$-differentiability is much less restrictive that the usual complex differentiability.
For example, $\L$-differentiability at $z_0$ does not imply continuity at $z_0$. However, $\L$-differentiability
in an open set $\Omega \subset \L$ implies usual differentiability in $\Omega$.
\end{remark}
Introducing the paracomplex operators:
\begin{equation}\label{dpc}
\frac{\partial}{\partial z} = \frac{1}{2}\Big(\frac{\partial}{\partial u} + \tau\frac{\partial}{\partial v}\Big),\qquad
\frac{\partial}{\partial \bar{z}} = \frac{1}{2}\Big(\frac{\partial}{\partial u} - \tau\frac{\partial}{\partial v}\Big),\nonumber
\end{equation}
where $z = u + \tau\, v$, we have that a differentiable function $f:\Omega\rightarrow\L$ is $\L$-differentiable if and only if
\begin{equation}\label{eqldif}
\displaystyle\frac{\partial f}{\partial \bar{z}} = 0.
\end{equation}
We observe that, writing  $f(u,v) = a(u,v) +\tau\, b(u,v)$, $u + \tau\,v\in\Omega$, the condition~\eqref{eqldif} is equivalent to the para-Cauchy-Riemann equations:
$$
\left\{\begin{aligned}
\frac{\partial a}{\partial u} &= \frac{\partial b}{\partial v},\\
\frac{\partial a}{\partial v} &= \frac{\partial b}{\partial u},
\end{aligned}\right.
$$
whose integrability conditions are given by the wave equations
$$a_{uu} - a_{vv} = 0 = b_{uu} - b_{vv}.$$

 We observe that the harmonicity condition is a system of partial differential equations (really an integral differential equation, since the $\Gamma$'s must be computed along a solution). Hence, in general, it is quite hard to find explicit solutions. However, for certain ambient spaces, as the Lie groups, these equations are essentially equivalent to a system of partial differential equations with constant coefficients. We will comment now, briefly, the case where the ambient space is a Lie group. In what follows $\K$ will denote either the complex numbers $\C$ or the Lorentz numbers $\L$.

Let $M$ be a $n$-dimensional Lie group endowed with a left-invariant Riemannian or Lorentzian metric  and let $f:\Omega\subset \K \rightarrow M$ be a conformal minimal immersion, where $\Omega \subset \K$ is an open set.
Let $\{e_1, e_2,\dots, e_n\}$ be a left-invariant orthonormal frame field, with $e_1,\dots,e_{n-1}$ spacelike and $e_n$ timelike if the metric is Lorentzian. We can write the (para)complex tangent field $\phi=\displaystyle\frac{\partial f}{\partial z}$ along $f$ both in terms of local coordinates $\{x_1,x_2,\dots,x_n\}$  in $M$ and, also, using the left-invariant vector fields. Hence, one has
\begin{equation}
\phi = \sum_{i=1}^{n} \phi_i\frac{\partial}{\partial x_i} = \sum_{i=1}^{n} \psi_i e_i,\nonumber
\end{equation}
where the functions $\phi_i$ and $\psi_i$\ are related by
\begin{equation}\label{eq:3.1}
\phi_i = \sum_{j=1}^{n} A_{ij}\psi_j,
\end{equation}
where $A: \Omega \rightarrow GL(n,\r)$ is a smooth map. In terms of the components $\psi_i$, the harmonicity condition can be written as
$$
\frac{\partial\psi_k}{\partial\bar{z}} + \frac{1}{2}\sum_{i,j=1}^{n} L_{ij}^k\, \bar{\psi_i}\,\psi_j = 0,
$$
where the symbols $L_{ij}^k$ are defined by
$$
\nabla_{e_i} e_j = \frac{1}{2}\sum_{k=1}^{n} L_{ij}^k\, e_k.
$$

Consequently, in the case of $n$-dimensional Lie groups, the Theorem~\ref{classicalwr} may be rephrased as follows
\begin{theorem}\label{teo3.6}
Let $M$ be a $n$-dimensional Lie group endowed with a left-invariant
Lorentzian metric and let $\{e_1,e_2,\dots,e_n\}$ be a left-invariant orthonormal frame field.
Let $f: \Omega \rightarrow M$ be a conformal minimal immersion, where $\Omega \subset\K$ is an open set.
We denote by $\phi \in \Gamma(f^\ast TM\otimes \K)$ the (para)complex tangent vector
$$
\phi = \frac{\partial f}{\partial z}.
$$
Then, the components $\psi_i$, $i= 1,\dots,n$, of $\phi$
satisfy the following conditions:
\begin{enumerate}\label{eq:3.2}
  \item [i)] $|\psi_1|^2 + |\psi_2|^2 + \cdots +|\psi_{n-1}|^2 - |\psi_n|^2 \neq 0$,
  \item [ii)]${\psi_1}^2 + {\psi_2}^2 + \cdots +\psi_{n-1}^2- {\psi_n}^2 = 0$,
  \item [iii)] $\displaystyle\frac{\partial\psi_k}{\partial\bar{z}} +
  \frac{1}{2}\sum_{i,j=1}^{n} L_{ij}^k\, \bar{\psi_i}\,\psi_j = 0$.
\end{enumerate}
Conversely, if $\Omega\subset\K$ is a simply connected domain and $\psi_k:\Omega\to\K$, $k = 1,\dots,n,$ are (para)complex functions satisfying the conditions above, then the map $f:\Omega\to M$ which coordinates are given by
$$
f_i = 2\,\Re\int  \sum_{j=1}^{n} A_{ij}\,\psi_j\, dz, \qquad i = 1,\dots,n,
$$
is a well-defined conformal minimal immersion.
\end{theorem}

\section{The Weierstrass representation in the Lorentzian Damek-Ricci spaces  $\S_4^1$}\label{quattro}

We consider the $4$-dimensional space $\S_4^1$ with global coordinates $\{x,y,z,t\}$. The left-invariant Lorentzian metric $g$ is given by:
$$
    g = e^{-t}\,dx^2 + e^{-t}\,dy^2 + e^{-2t}\,(dz + \frac{c}{2}y\,dx - \frac{c}{2}x\,dy)^2 - \,dt^2,
$$
where $c \in \r$. The Lie algebra $\s_4$ of $\S_4^1$ has an orthonormal basis
$$
 e_1 = e^{\frac{t}{2}}\Big{(} \frac{\partial}{\partial x} - \frac{c\,y}{2}\frac{\partial}{\partial z}\Big{)}, \qquad
 e_2 = e^{\frac{t}{2}}\Big{(} \frac{\partial}{\partial y} + \frac{c\,x}{2}\frac{\partial}{\partial z}\Big{)},\qquad
 e_3 = e^{t}\frac{\partial}{\partial z}, \qquad
 e_4 = \frac{\partial}{\partial t},
$$
where $e_1,e_2,e_3$ are spacelike and $e_4$ is timelike.
The Lie brackets are given by
\begin{equation}
\left\{
\begin{aligned}
\displaystyle [e_1,e_2] &= c\,e_3,\quad [e_1,e_3] = 0,\quad [e_1,e_4] = -\frac{1}{2}e_1, \\
\displaystyle [e_2,e_3] &= 0, \quad [e_2,e_4] = -\frac{1}{2}e_2, \quad [e_3,e_4] = -e_3.
\end{aligned}
\right.
\nonumber
\end{equation}
The Levi-Civita connection is given by:
$$\begin{aligned}
  \nabla_{e_1} e_1 &= -\dfrac{1}{2}e_4,\quad \nabla_{e_1} e_2 = \dfrac{c}{2}e_3,\quad
  \nabla_{e_1} e_3 = -\dfrac{c}{2}e_2,\quad \nabla_{e_1} e_4 = -\dfrac{1}{2}e_1,\\
  \nabla_{e_2} e_1 &= -\dfrac{c}{2}e_3,\quad \nabla_{e_2} e_2 = -\dfrac{1}{2}e_4,\quad \nabla_{e_2} e_3 = \dfrac{c}{2}e_1,\quad \nabla_{e_2} e_4 = -\dfrac{1}{2}e_2,
 \\ \nabla_{e_3} e_1 &= -\dfrac{c}{2}e_2,\quad \nabla_{e_3} e_2 = \dfrac{c}{2}e_1,\quad
\nabla_{e_3} e_3 = - e_4,\quad \nabla_{e_3} e_4 = -e_3,
 \\ \nabla_{e_4} e_1& = \nabla_{e_4} e_2 = \nabla_{e_4} e_3 = \nabla_{e_4} e_4 = 0.
\end{aligned}
$$
Also, we have that the non zero $L_{ij}^k$ are:
$$
\begin{aligned}
\displaystyle L_{11}^4 &= -1,\quad L_{12}^3 = c, \quad L_{13}^2 = -c, \quad L_{14}^1 = -1,
\quad L_{21}^3 = -c, \quad L_{22}^4 = -1,\\
\displaystyle L_{23}^1 & = c,\quad L_{24}^2 = -1, \quad L_{31}^2 = -c, \quad L_{32}^1 = c,
\quad L_{33}^4 = -2, \quad L_{34}^3 = -2.
\end{aligned}
$$
The matrix $A$ defined in the previous section is
\begin{equation}
A = \left[
  \begin{array}{cccc}
            e^{\frac{t}{2}}      &           0                 &  0  & 0 \\
                 0               &     e^{\frac{t}{2}}         &  0  & 0 \\
    -\frac{c}{2}\,e^{\frac{t}{2}}y & \frac{c}{2}\,e^{\frac{t}{2}}x & e^t & 0 \\
                 0               &            0                &  0  & 1
  \end{array}
\right].\nonumber
\end{equation}
The harmonicity condition is given by the following system of PDEs:
\begin{equation}\label{eq:4.11} \left\{\begin{aligned}
\displaystyle\frac{\partial\psi_1}{\partial\bar{z}}& -\frac{1}{2}\bar{\psi_1}\psi_4 + c\,\Re(\bar{\psi_2}\psi_3) = 0,\\
\displaystyle\frac{\partial\psi_2}{\partial\bar{z}}& -\frac{1}{2}\bar{\psi_2}\psi_4 - c\,\Re(\bar{\psi_1}\psi_3) = 0,\\
\displaystyle\frac{\partial\psi_3}{\partial\bar{z}}& -\bar{\psi_3}\psi_4 + \frac{c}{2}\,(\bar{\psi_1}\psi_2 -
\bar{\psi_2}\psi_1) = 0,\\
\displaystyle\frac{\partial\psi_4}{\partial\bar{z}}& - \frac{1}{2}(\bar{\psi_1}\psi_1 + \bar{\psi_2}\psi_2) - \bar{\psi_3}\psi_3 = 0.
\end{aligned}\right.
\end{equation}
Then Theorem~\ref{teo3.6} takes the form
\begin{theorem}\label{TdaRdeWpS4}
Let $\Omega \subseteq \K$ be an open set endowed of a (para)complex coordinates
and $f: \Omega\rightarrow\S_4^1$ a conformal minimal immersion. Then, the components of the (para)complex tangent vector
$$
\phi = \frac{\partial f}{\partial z} =
\sum_{i=1}^{4} \psi_i\, e_i$$
satisfy the system \eqref{eq:4.11} and the following conditions:
\begin{enumerate}
  \item [i)] $|\psi_1|^2 + |\psi_2|^2 + |\psi_3|^2 - |\psi_4|^2 \neq 0$,
  \item [ii)]${\psi_1}^2 + {\psi_2}^2 + {\psi_3}^2 - {\psi_4}^2= 0$.
\end{enumerate}
Conversely, if $\Omega$ is simply connected and $\psi_i: \Omega \rightarrow \K$, $i=1, 2, 3,4$, are functions
satisfying the above conditions, then the map $f: \Omega\rightarrow\S_4^1$ with coordinates
$$
f_i = 2\,\Re\int  \sum_{j=1}^{4} A_{ij}\psi_j\, dz,\qquad i = 1,2,3,4,
$$
defines a conformal minimal immersion in $\S_4^1$.
\end{theorem}
We will give now some examples.

\begin{example}
Consider the paracomplex functions
$$\psi_1 = \frac{\tau}{u},\qquad \psi_2 = \psi_3 = 0, \qquad\psi_4 = \frac{1}{u},$$
defined in the simply connected domain $\Omega =\{u + \tau\,v \in \L \;|\; u> 0\}$.
We have that \eqref{eq:4.11} and conditions i) and ii) of Theorem~\ref{TdaRdeWpS4} are satisfied and, so,
the map $f:\Omega \rightarrow \S_4^1$ given by
$$\left\{\begin{aligned}
\displaystyle f_1 &= \frac{2(v - v_0)}{u_0},\\
\displaystyle f_2 &= k\in\r,\\
\displaystyle f_3 &= -\frac{k_1(v - v_0)}{u_0}, \quad k_1\in\r,\\
\displaystyle f_4 &= 2\ln{\bigg{(}\frac{u}{u_0}\bigg{)}},
\end{aligned}\right.
$$
is a conformal timelike minimal immersion, where $z_0 = u_0 + \tau\,v_0 \in \Omega$.
\end{example}
\begin{example}
The paracomplex functions
$$\psi_1 = \psi_2 = \frac{\tau}{\sqrt{2}u},\quad\psi_3 = 0,\qquad \psi_4 = \frac{1}{u}$$ defined in $\Omega =\{u + \tau\,v \in \L\;|\; u> 0\}$
satisfy the Theorem~\ref{TdaRdeWpS4}. Then, the map $f:\Omega \rightarrow \S_4^1$ with coordinates
$$\left\{\begin{aligned}
\displaystyle f_1 &= -\frac{2(v - v_0)}{u_0},\\
\displaystyle f_2 &= -\frac{2(v - v_0)}{u_0},\\
\displaystyle f_3 &= k\in\r,\\
\displaystyle f_4 &= 2\ln{\bigg{(}\frac{u}{u_0}\bigg{)}},
\end{aligned}\right.$$
is a conformal timelike minimal immersion, where $z_0 = u_0 + \tau\,v_0 \in \Omega$.
\end{example}

\section{The Weierstrass representation in the Lorentzian Damek-Ricci spaces $\S_4^3$}\label{cinque}

We consider the $4$-dimensional Lorentzian Damek-Ricci space $\S_4^3$ with global coordinates $\{x,y,z,t\}$. The left-invariant Lorentzian metric $g$ on $\S_4^3$ is given by:
\begin{equation}
    g = e^{-t}\,dx^2 + e^{-t}\,dy^2 - e^{-2t}\,\Big(dz + \frac{c}{2}y\,dx - \frac{c}{2}x\,dy\Big)^2 + dt^2,\nonumber
\end{equation}
where $c \in \r$. The Lie algebra $\s_4$ of $\S_4^3$ has the orthonormal basis
$$
 e_1 = e^{\frac{t}{2}}\big{(} \frac{\partial}{\partial x} - \frac{c\,y}{2}\frac{\partial}{\partial z}\big{)}, \qquad
 e_2 = e^{\frac{t}{2}}\big{(} \frac{\partial}{\partial y} + \frac{c\,x}{2}\frac{\partial}{\partial z}\big{)}, \qquad
 e_3 = e^{t}\frac{\partial}{\partial z}, \qquad
 e_4 = \frac{\partial}{\partial t},$$
where $e_1,e_2,e_4$ are spacelike and $e_3$ is timelike.
The Lie brackets are given by
\begin{equation}
\begin{aligned}\nonumber
[e_1,e_2] &=c\,e_3, \qquad [e_1,e_3] = 0,\quad [e_1,e_4] = -\frac{1}{2}e_1, \\
[e_2,e_3] &= 0, \qquad [e_2,e_4] = -\frac{1}{2}e_2, \quad [e_3,e_4] = -e_3.
\end{aligned}
\end{equation}
As the Levi-Civita connection is given by:
$$
\begin{aligned}
   \nabla_{e_1} e_1 &= \dfrac{1}{2}\,e_4,\qquad \nabla_{e_1} e_2 = \dfrac{c}{2}\,e_3,\qquad
  \nabla_{e_1} e_3 = \dfrac{c}{2}\,e_2,\qquad \nabla_{e_1} e_4 = -\dfrac{1}{2}\,e_1,
 \\ \nabla_{e_2} e_1& = -\dfrac{c}{2}\,e_3,\qquad \nabla_{e_2} e_2 = \dfrac{1}{2}\,e_4,\qquad \nabla_{e_2} e_3 = -\dfrac{c}{2}\,e_1,\qquad \nabla_{e_2} e_4 = -\dfrac{1}{2}\,e_2,
 \\ \nabla_{e_3} e_1 &= \dfrac{c}{2}\,e_2,\qquad \nabla_{e_3} e_2 = -\dfrac{c}{2}e_1,\qquad
\nabla_{e_3} e_3 = -e_4,\qquad \nabla_{e_3} e_4 = -e_3,
 \\ \nabla_{e_4} e_1 &= \nabla_{e_4} e_2 = \nabla_{e_4} e_3 = \nabla_{e_4} e_4 = 0,
\end{aligned}
$$
then the non zero $L_{ij}^k$ are
$$
\begin{aligned}
\displaystyle L_{11}^4 &= 1,\quad L_{12}^3 = c, \quad L_{13}^2 = c, \quad L_{14}^1 = -1,
\quad L_{21}^3 = -c, \quad L_{22}^4 = 1, \\
\displaystyle L_{23}^1 & = -c,\quad L_{24}^2 = -1, \quad L_{31}^2 = c, \quad L_{32}^1 = -c,
\quad L_{33}^4 = -2, \quad L_{34}^3 = -2.
\end{aligned}
$$
Also the matrix $A$ is given by
$$
A = \left[
  \begin{array}{cccc}
            e^{\frac{t}{2}}      &           0                 &  0  & 0 \\
                 0               &     e^{\frac{t}{2}}         &  0  & 0 \\
    -\frac{c}{2}\,e^{\frac{t}{2}}y & \frac{c}{2}\,e^{\frac{t}{2}}x & e^t & 0 \\
                 0               &            0                &  0  & 1
  \end{array}
\right].$$
The harmonicity condition becomes
\begin{equation}\label{eq:4.12} \left\{\begin{aligned}
\displaystyle\frac{\partial\psi_1}{\partial\bar{z}}& -\frac{1}{2}\bar{\psi_1}\psi_4 - c\,\Re(\bar{\psi_2}\psi_3) = 0,\\
\displaystyle\frac{\partial\psi_2}{\partial\bar{z}}& -\frac{1}{2}\bar{\psi_2}\psi_4 + c\,\Re(\bar{\psi_1}\psi_3) = 0,\\
\displaystyle\frac{\partial\psi_3}{\partial\bar{z}}& -\bar{\psi_3}\psi_4 + \frac{c}{2}\,(\bar{\psi_1}\psi_2 -
\bar{\psi_2}\psi_1) = 0,\\
\displaystyle\frac{\partial\psi_4}{\partial\bar{z}}& + \frac{1}{2}(\bar{\psi_1}\psi_1 + \bar{\psi_2}\psi_2) -\bar{\psi_3}\psi_3 = 0.
\end{aligned}\right.
\end{equation}
Therefore, the Weierstrass representation formula is given by the following
\begin{theorem}\label{TdaRdeWpS1}
Let $\Omega \subseteq \K$ be an open set endowed of a (para)complex coordinate
and $f: \Omega \rightarrow \S_4^3$ a conformal minimal immersion. Then, the components of the (para)complex tangent vector
$$
\phi = \frac{\partial f}{\partial z} =
\sum_{i=1}^{4} \psi_i\, e_i,$$
satisfy the system \eqref{eq:4.12} and the following conditions:
\begin{enumerate}
  \item [i)] $|\psi_1|^2 + |\psi_2|^2 - |\psi_3|^2 + |\psi_4|^2 \neq 0$,
  \item [ii)]${\psi_1}^2 + {\psi_2}^2 - {\psi_3}^2 + {\psi_4}^2= 0$.
\end{enumerate}

Conversely, if $\Omega$ is simply connected and $\psi_i: \Omega \rightarrow \K, i=1, 2, 3,4$, are functions
satisfying the above conditions, then the map $f: \Omega \rightarrow \S_4^3$ with coordinates
$$
f_i = 2\,\Re\int  \sum_{j=1}^{4} A_{ij}\,\psi_j\, dz,\qquad i = 1,2,3,4,
$$
is a conformal timelike (or spacelike) minimal immersion in $\S_4^3$.
\end{theorem}

We will give now some examples.
\begin{example}
It easy to check that the complex functions
$$\psi_1 = \frac{i}{u},\quad \psi_2 = \psi_3 = 0, \qquad\psi_4 = \frac{1}{u},$$
defined in $\Omega =\{u + i\,v \in \C \;|\; u> 0\}$,
satisfy \eqref{eq:4.12} and the conditions i) and ii) of Theorem~\ref{TdaRdeWpS1}. So, the map $f:\Omega \rightarrow \S_4^3$ given by:
$$\left\{\begin{aligned}
\displaystyle f_1 &= -\frac{2(v - v_0)}{u_0},\\
\displaystyle f_2 &= k\in\r,\\
\displaystyle f_3 &= \frac{k_1(v - v_0)}{u_0},\quad k_1\in\r,\\
\displaystyle f_4 &= 2\ln{\bigg{(}\frac{u}{u_0}\bigg{)}},
\end{aligned}\right.$$
is a conformal spacelike minimal immersion, where $z_0 = u_0 + i\,v_0 \in \Omega$.
\end{example}

\begin{example}
Consider the paracomplex functions
$$\psi_1 = \psi_2 = 0,\quad\psi_3 = \frac{\tau}{2u},\qquad\psi_4 = \frac{1}{2u}$$
in the simply connected domain $\Omega =\{u + \tau\,v \in \L \;|\; u> 0\}$.

We have that they satisfy \eqref{eq:4.12} and the conditions i) and ii) of Theorem~\ref{TdaRdeWpS1}. Then, the map $f:\Omega \rightarrow \S_4^3$ with components:
$$\left\{\begin{aligned}
\displaystyle f_1 &= k_1\in\r,\\
\displaystyle f_2 &= k_2\in\r,\\
\displaystyle f_3 &= \frac{(v - v_0)}{u_0},\\
\displaystyle f_4 &= \ln{\bigg{(}\frac{u}{u_0}\bigg{)}},
\end{aligned}\right.$$
is a conformal timelike minimal immersion in $\S_4^3$, where $z_0 = u_0 + \tau\,v_0 \in \Omega$.
\end{example}

\end{document}